\documentclass[a4paper,11pt]{amsart}
\pdfoutput=1
\usepackage{graphicx}
\usepackage[ansinew]{inputenc}
\usepackage[T1]{fontenc}
\usepackage{amssymb}
\usepackage{mathrsfs}
\usepackage{amsfonts,amscd,amsmath}
\usepackage{amssymb,latexsym}
\usepackage[all]{xy}
\usepackage{tikz-cd} 

\newtheorem{thm}{Theorem}[section]

\newtheorem{lem}[thm]{Lemma}
\newtheorem{prop}[thm]{Proposition}

\theoremstyle{definition}
\newtheorem{definition}[thm]{Definition}

\newtheorem{rem}[thm]{Remark}
\newtheorem{ex}[thm]{Example}
\usepackage{mystyle}

\usepackage{color}

\usepackage{hyperref}
\hypersetup{colorlinks,%
citecolor=black,%
filecolor=black,%
linkcolor=black,%
urlcolor=black,%
pdftex}

\def\pd{\operatorname{pd}}

\begin{document}

\author{Juan J. Nu\~no-Ballesteros, Guillermo Pe\~nafort Sanchis, Cinzia Villa}

\address{Departament de Matem\`atiques,
Universitat de Val\`encia, Campus de Burjassot, 46100 Burjassot,
Spain}
\email{Juan.Nuno@uv.es}
\email{Guillermo.Penafort@uv.es}
\email{cinvi@alumni.uv.es}
\title[Target double and triple points of multi-germs]{On the analytic structure of double  and  triple points in the target of finite holomorphic multi-germs}

%

\begin{abstract}
We study the analytic structure of the double and triple point spaces $M_2(f)$ and $M_3(f)$ of finite multi-germs $f\colon (X,S)\to(\C^{n+1},0)$, based on results of Mond and Pellikaan for the mono-germ case. We show that these spaces are Cohen-Macaulay, provided that certain dimensional conditions are satisfied, and give explicit expressions for their defining ideals in terms of those of their mono-germ branches.
  \end{abstract}

								\maketitle
								\section{Introduction}
In this text we study double and triple points of finite holomorphic multi-germs 
\[f\colon (X,S)\to (\C^{n+1},0),\]
where $X$ stands for a Cohen-Macaulay analytic space of dimension $n$. By saying that $(X,S)$ is a multi-germ, the set $S$ is assumed to be finite, and by saying  that $X$ is an analytic space (in contrast to a complex space), it is understood that $X$ is reduced. These assumptions and the fact that $f$ is holomorphic will not be repeated every time these objects appear, but  extra hypothesis will sometimes be added on top of them.

There are several ways to define double, triple and, in general, $k$-multiple point spaces, and this work focusses on the spaces $M_k(f)\subseteq (\C^{n+1},0)$, defined by means of Fitting ideals of the pushforward module $f_*\cO_{(X,S)}$. These spaces were introduced by Mond and Pellikaan \cite{Mond-Pellikaan:1989} , who studied their analytic structures in the mono-germ case. They show that, for mono-germs, the double and triple point spaces $M_2(f)$ and $M_3(f)$ are determinantal if they have the right dimension. As it turns out  their proofs do not apply in the multi-germ case, as explained in Example \ref{exNoSubmatrix}. The goal of this work is to study the analytic structure of $M_2(f)$ and $M_3(f)$ for multi-germs, giving conditions for them to be Cohen-Macaulay and showing explicit expressions for their ideals in terms of those of the spaces $M_k$ of the mono-germ branches of $f$.

 For a finite multi-germ $f$ as above, the pushforward $f_*\cO_{(X,S)}$ is a finitely presented $\cO_{n+1}$-module, and it admits a square presentation matrix $\lambda$ \cite[Section 2.2]{Mond-Pellikaan:1989}. Assuming $\lambda$ to be of size $q\times q$, the $k$th Fitting ideal of $f_*\cO_{(X,S)}$ is defined as follows:
\[F_k(f_*\cO_{(X,S)})=
\begin{cases}
 \langle \text{$q-k$ size minors of $\lambda$}\rangle,&\text{ if $k<q$}\\
 \cO_{n+1},&\text{otherwise}
\end{cases}
\]
Then, the  \emph{(target) $k$-multiple point space} of $f$  is the complex subspace $M_k(f)\subset(\C^{n+1},0)$, defined as
\[M_k(f)=V(F_{k-1}(f_*\cO_{(X,S)})).\]

This definition is based on the fact that, set-theoretically, $M_k(f)$ is the germ of the set of points $y\in \C^{n+1}$ where $(f_*\cO_X)_y$ requires at least a number $k$ of generators \cite[Proposition 1.5]{Mond-Pellikaan:1989}. Geometrically, this means that $M_k(f)$ is given by the points $y\in \C^{n+1}$ with at least $k$ preimages, counting with multiplicity. Observe that $M_1(f)$ is the hypersurface with equation $\det\lambda$ in $(\C^{n+1},0)$ and  coincides, as a set-germ, with the image $(f(X),0)$.

\begin{definition}
Given $f\colon (X,S)\to (\C^{n+1},0)$, we say that $M_k(f)$ is \emph{dimensionally correct} if $\dim M_k(f)=n-k+1$.
\end{definition}

From the fact that $(X,S)$ is assumed to be reduced and $f$ is assumed to be finite, it follows that $f$ must be locally an embedding on an open dense subset of $X$. From this it follows that $M_2(f)$ is dimensionally correct if and only if $f$ is generically one-to-one. Moreover,
$f$ is generically one-to-one if and only if $M_1(f)$ is reduced \cite[Proposition 3.2]{Mond-Pellikaan:1989}.
	
\begin{thm}[Mond-Pellikaan]\label{MondPellikaanDouble}
Any mono-germ $f\colon (X,0)\to(\C^{n+1},0)$ satisfies $\dim M_2(f)\geq n-1$. If the equality holds, then $M_2(f)$ is Cohen-Macaulay. If  furthermore $M_3(f)=\emptyset$,  then $M_2(f)$ is a complete intersection.
\end{thm}

The proof of this result consists on showing that the presentation matrix can be chosen to be symmetric and in such a way that  $F_1(f_*\cO_{X,0})$ is generated by minors of the submatrix obtained by removing the first row  \cite[Theorem 3.4]{Mond-Pellikaan:1989}. This shows $M_2(f)$ to be determinantal whenever it has dimension $n-1$ (observe that in codimension two the determinantal and Cohen-Macaulay conditions are equivalent, by the Hilbert-Burch theorem). The complete intersection case also follows from the fact that we can remove the first row, taking into account that the hypothesis that $M_3(f)=\emptyset$ is equivalent to the fact that the above matrix has size $2\times 2$.

In the case where $X$ is Gorenstein, the presentation matrix may be chosen in such a way that $F_2(f_*\cO_{X,0})$ is generated by minors of the submatrix obtained by removing the first row and first column \cite[Theorem 4.3]{Mond-Pellikaan:1989}. This implies what follows:

\begin{thm}[Mond-Pellikaan]\label{MondPellikaanTriple}
Any mono-germ $f\colon (X,0)\to(\C^{n+1},0)$, with $(X,0)$ Gorenstein, satisfies $\dim M_3(f)\geq n-2$. If the equality holds, then $M_3(f)$ is Cohen-Macaulay.
\end{thm}

Unfortunately, in the case of multi-germs, sometimes no proper submatrix of the presentation matrix has the right  minors to generate 
$F_1(f_*\cO_{(X,S)})$ nor $F_2(f_*\cO_{(X,S)})$. This is illustrated by the following example:

\begin{ex}\label{exNoSubmatrix}Consider the multi-germ $f\colon \bigsqcup_{i=1}^3(\C^2,0)\to (\C^3,0)$, parametrizing the planes $\{X=0\}$, $\{Y=0\}$ and $\{Z=0\}$ in $(\C^3,0)$. Since the branches are immersive mono-germs, their presentation matrices are $1\times1$ matrices whose only entry are the generators of the ideals of the corresponding planes. The presentation matrix of $f_*\cO_{(X,S)}$ is then the diagonal matrix 
	\[
	\lambda=
	\begin{bmatrix}
		X &0&0  \\
		0 &Y&0\\
		0& 0 & Z  \\
	\end{bmatrix}.
\]
As expected, $M_2(f)$ is the union of the coordinate axes of $(\C^3,0)$, and $M_3(f)$ is the origin, both spaces with reduced structure. One sees that no proper submatrix of $\lambda$ has the adequate entries to generate the corresponding equations.
\end{ex}

\begin{rem}
There are multiple point spaces $M_k(f)$ for  finite (global) mappings $f\colon X\to Y$, defined by the corresponding sheaf theoretical versions $\mathcal F_k$ of the Fitting ideals $F_k$. The properties that we show for  multi-germs imply the corresponding properties for mappings, because the stalks of $\mathcal F_k$ are the ideals $F_k$ associated to a multi-germ. We omit these translations to the global setting, as they are quite obvious.
\end{rem}

								\section{Target multiple-point spaces for multi-germs}

To simplify writing, for the remainder of the text we fix some notations regarding multi-germ branches. Given $f\colon (X,S)\to (\C^{n+1},0)$, with $S=\{x^{(1)},\dots,x^{(r)}\}$, for each $i\in \{1,\dots,r\}$, we write  \[f^{(i)}=f\vert_{(X,x^{(i)})},\qquad F_j^{(i)}=F_j(f^{(i)}_*\cO_{X,x^{(i)}})\quad \text{and}\quad M_j^{(i)}=M_j(f^{(i)}).\]

From the equality $f_*\cO_{(X,S)}=\bigoplus_{i=1}^r f^{(i)}_*\cO_{X,x^{(i)}}$, one observes that the associated Fitting ideals are generated by minors of a block diagonal matrix, and then
\[F_k(f_*\cO_{(X,S)})=\sum_{j_1+\dots+j_r= k}F_{j_1}^{(1)}\cdot \ldots\cdot F_{j_r}^{(r)}.\]
In combination with the inclusions $M_j^{(i)}\subseteq M_{j-1}^{(i)}$, the previous equality gives the \emph{set theoretical} decomposition
\begin{equation}\label{eqMkDecomposition}
M_k(f)=\bigcup_{\substack{\smallskip1\leq \ell\leq r\\\smallskip1\leq i_1<\dots <i_\ell\leq r\\ j_1+\dots+j_\ell=k}} M_{j_1}^{(i_1)}\cap\ \dots\  \cap M_{j_\ell}^{(i_\ell)}.
\end{equation}
This decomposition corresponds to the intuitive idea of what multiple points of multi-germs should be. For example, it states that, as a set, the triple points of a bigerm are the union of the triple points of the first branch, the intersection of the double points of the first branch with the image of the second, the intersection of the image of the first branch with the double points of the second and, finally, the triple points of the second branch. 

It is natural to ask to what extent this decomposition gives information about $M_k(f)$ as a complex space. To be precise, we want to know if the ideal $F_k(f_*\cO_{(X,S)})$ can be expressed as
\[F_k(f_*\cO_{(X,S)})=\bigcap_{\substack{\smallskip1\leq \ell\leq r\\\smallskip1\leq i_1<\dots <i_\ell\leq r\\ j_1+\dots+j_\ell=k-\ell+1}} F_{j_1}^{(i_1)}+\dots+F_{j_\ell}^{(i_\ell)}.\]

We will use the following three basic lemmas. The first follows easily from localizing at associated primes and the second follows from the interpretation of $\depth$ in terms of $\Ext$, by looking at the long exact sequence associated to a short exact sequence (see e.g. \cite[\href{https://stacks.math.columbia.edu/tag/00LX}{Tag 00LX}]{stacks-project}). For completeness, we include a proof of the third lemma. 							

\begin{lem}\label{lemSumInterCommute}
Let $R$ be a noetherian local ring, and let $J_1,\dots,J_r$ and $I$ be ideals in $R$. Assume that for all $i\neq j$ the condition $J_i\not\subseteq \kp$ holds for all $\kp \in \Ass(J_j+I)$. Then, $\bigcap_i (J_i+I)=(\bigcap_iJ_i)+I$.
\end{lem}

\begin{lem}[Depth Lemma]\label{lemDepthLemma}
Let $R$ be a local Noetherian ring and let $0\to N'\to N\to N''\to 0$ be a short exact sequence of nonzero finitely generated $R$-modules.
\begin{enumerate}
\item \label{lemDepthLemmaItem1} $\depth N\geq \min\{\depth N', \depth N''\}$.
\item \label{lemDepthLemmaItem2} $\depth N''\geq \min\{\depth N, \depth N'-1\}$.
\item \label{lemDepthLemmaItem3} $\depth N'\geq \min\{\depth N, \depth N''+1\}$.
\end{enumerate}
\end{lem}

\begin{lem}\label{lemProductIdealsDepth} Let $I$ and $J$ be ideals in a local noetherian ring $R$, such that $\depth\frac{R}{I}< \depth R$, $\depth\frac{R}{J}< \depth R$ and $\depth\frac{R}{IJ}< \depth R$. Then 
\[\depth \frac{R}{IJ}\geq \depth\frac{R}{I}+\depth\frac{R}{J}-\depth R+1.\]
\begin{proof}
The hypotheses imply $\depth  \frac{R}{I}=\depth I-1$, $\depth  \frac{R}{J}=\depth J-1$ and $\depth  \frac{R}{IJ}=\depth IJ-1$. By the third of these equalities, it suffices to show that $\depth IJ\geq \depth \frac{R}{I}+\depth \frac{R}{J}-\depth R+2$ and, by the Auslander-Buchsbaum formula, this amounts to showing $\pd IJ\leq 2\depth R-\depth \frac{R}{I}-\depth \frac{R}{J}-2$. Again by the Auslander-Buchsbaum formula, we obtain
\begin{align*}
\pd IJ=\pd (I\otimes J)&\leq\pd I+\pd J=2\depth R-\depth I-\depth J\\&=2\depth R-\depth \frac{R}{I}-\depth \frac{R}{J}-2.\qedhere
\end{align*}
\end{proof}
\end{lem}

Now we are ready to prove that, when they are dimensionally correct, double points of multi-germs are Cohen-Macaulay and to show how to compute them under favorable circumstances:

\begin{thm}\label{thmMultigermDoublePoints}If $f\colon (X,S)\to (\C^{n+1},0)$ is finite and generically one-to-one, then $M_2(f)$ is Cohen-Macaulay of dimension $n-1$. If furthermore $M_3(f)$ is dimensionally correct, then the ideal of $M_2(f)$ satisfies
\[F_1(f_*\cO_{(X,S)})=\bigcap_{1\leq i\leq r}F_1^{(i)}\cap \bigcap_{1\leq i<j\leq r}\big(F_0^{(i)}+F_0^{(j)}\big).\]

\begin{proof}	We proceed by induction over $r=|S|$. If $r=1$, then $M_2(f)$ is Cohen-Macaulay by Theorem \ref{MondPellikaanDouble} and the expression for the defining ideal reads $F_1(f_*\cO_{(X,S)})=F_1(f_*\cO_{(X,S)})$, so there is nothing to prove. Now assume the statement to be true for multi-germs with no more than $r-1$ branches.

Showing $M_2(f)$ to be Cohen-Macaulay amounts to showing \[\depth \frac{\cO_{n+1}}{F_1(f_*\cO_{X,S})}=n-1,\] because the generically one-to-one condition implies $\dim M_2(f)=n-1$. To exploit the inductive hypothesis, consider $S'=\{x^{(1)},\dots,x^{(r-1)}\}$ and let \[f'= f_{(X,S')}\quad\text{and}\quad F'_j=F_j(f'_*\cO_{(X,S')}),\] so that $F_1(f_*\cO_{X,S})=F_1' F_0^{(r)}+F_0' F_1^{(r)}$. 

Since $F_0'$ and $F_0^{(r)}$ are principal, it follows from Lemma \ref{lemProductIdealsDepth} and the induction hypothesis that \[\depth \frac{\cO_{n+1}}{F_1'  F_0^{(r)}}\geq n-1\quad\text{and}\quad\depth \frac{\cO_{n+1}}{F_0'  F_1^{(r)}}\geq n-1.\]  Applying Item (\ref{lemDepthLemmaItem2}) of the \hyperref[{lemDepthLemma}]{Depth Lemma} to the short exact sequence
\[0\to \frac{\cO_{n+1}}{F_1'  F_0^{(r)}\cap F_0'  F_1^{(r)}}\to \frac{\cO_{n+1}}{F_1'  F_0^{(r)}}\oplus \frac{\cO_{n+1}}{F_0'  F_1^{(r)}}\to \frac{\cO_{n+1}}{F_1'  F_0^{(r)}+ F_0'  F_1^{(r)}}\to 0\]
reduces the problem to showing that \[\depth\frac{\cO_{n+1}}{F_1'  F_0^{(r)}\cap F_0'  F_1^{(r)}}\geq n.\]
From the inclusions $F_0'\subseteq F_1'$ and $F_0^{(r)}\subseteq F_1^{(r)}$, one obtains
	\[F'_1 F ^{(r)}_0 \cap  F'_0  F ^{(r)}_1\subseteq  F'_1\cap  F ^{(r)}_0\cap  F'_0\cap  F ^{(r)}_1= F'_0\cap  F ^{(r)}_0,\]
	\[F'_{1} F ^{(r)}_{0} \cap F'_{0}  F ^{(r)}_{1} \supseteq F'_0 F ^{(r)}_0.\] 
Since $F_0'$ and $F_0^{(r)}$ are principal ideals defining the images of $f'$ and $f^{(r)}$, the generically one-to-one condition forces $F'_0 F ^{(r)}_0=F'_0\cap F ^{(r)}_0$. From the two inclusions above, now we obtain $F'_1 F ^{(r)}_0 \cap F'_0  F ^{(r)}_1=F'_0 F ^{(r)}_0$. Being the ideal of a hypersurface in $(\C^{n+1},0)$, its quotient must have depth $n$, as desired. This finishes the proof that $M_2(f)$ is Cohen-Macaulay.
 
 Now assume $M_3(f)$ to be dimensionally correct. This forces $F_0'$ and $F_0^{(r)}$ to be generated by regular elements in $\cO_{n+1}/F_1^{(r)}$ and $\cO_{n+1}/F'_1$. This is true because $M_2(f')$ and $M_2^{(r)}$ are Cohen-Macaulay of dimension $n-1$, and points in $M_2(f')\cap M_1(f^{(r)})$ or $M_1(f')\cap M_2(f^{(r)})$ are triple points of $f$, which are assumed to be have dimension $(n-2)$.
 Then, we obtain
 \begin{align*}
 F_1(f_*\cO_{(X,S)})	&=F_1'F_0^{(r)}+F_0'F_1^{(r)}\\
 				&=\big(F_1'\cap F_0^{(r)}\big)+\big(F_0'\cap F_1^{(r)}\big)\\
				&=F_1'\cap F_1^{(r)}\cap \big(F_0'+F_0^{(r)}\big)\\
				&=\bigcap_{1\leq i\leq r-1}F_1^{(i)}\cap \bigcap_{1\leq i<j\leq r-1}\big(F_0^{(i)}+F_0^{(j)}\big)\cap F_1^{(r)}\cap \big(F_0'+F_0^{(r)}\big),
 \end{align*}
where the third equality follows from the inclusions $F_0'\subseteq F_1'$ and $F_0^{(r)}\subseteq F_1^{(r)}$ and the last equality follows from the induction hypothesis (observe that $M_3(f')$ is dimensionally correct, because every ingredient in decomposition of $M_3(f')$ in (\ref{eqMkDecomposition}) appears in the decomposition of $M_3(f)$ as well). Comparign this last expression to the formula in the statement shows that it is enough to prove that \[F_0'+F_0^{(r)}=\bigcap_{1\leq i\leq r-1}\big(F_0^{(i)}+F_0^{(r)}\big).\]
The generically one-to-one condition implies \[F_0'=\bigcap_{1\leq i\leq r-1} F_0^{(i)}\] and the fact that $F_0^{(r)}$ is generated by a regular element in $\cO_{n+1}/F_0'$. The desired equality now follows, taking into account that the ideals $F_0^{(i)}$ are principal.
\end{proof}
\end{thm}

\begin{ex}Let  $f\colon (\C^1,0)\sqcup(\C^1,0)\sqcup (\C^1,0)\to (\C^2,0)$ be the map-germ with branches $f^{(1)}(x)=(x,0)$, $f^{(2)}(x)=(0,x)$ and $f^{(3)}(x)=(x,x)$. Since $f^{(i)}$ is immersive, it follows that $F_0^{(i)}$ is the reduced ideal of the image of $f^{(i)}$ and $F_1^{(i)}=\cO_2$.
Therefore,
\begin{align*}
F_1(f_*(\cO_{(X,S)})&=F_1^{(1)}F_0^{(2)}F_0^{(3)}+F_0^{(1)}F_1^{(2)}F_0^{(3)}+F_0^{(1)}F_0^{(2)}F_1^{(3)}\\
				&=\langle X(X-Y), Y(X-Y), XY\rangle.
\end{align*}
At the same time, the right hand side of the equality in Theorem \ref{thmMultigermDoublePoints} gives
\[F_1^{(1)}\cap F_1^{(2)}\cap F_1^{(2)}\cap \big(F_0^{(1)}+F_0^{(2)}\big)\cap \big(F_0^{(1)}+F_0^{(3)}\big)\cap \big(F_0^{(2)}+F_0^{(3)}\big)=\langle X,Y\rangle.\]
The disagreement comes from the fact that the origin is a triple point of $f$ and, in these dimensions, $M_3(f)$ must be empty to be dimensionally correct.
\end{ex}

To study triple points, we will use a typical trick in the theory of singular mappings, consisting on replacing the original mapping by an unfolding having desired properties. When studying map-germs $f\colon (\C^n,S)\to (\C^{n+1},0)$, an unfolding of $f$ is simply a map germ of the form $(f_t,t)\colon  (\C^n,S)\times (\C^r,0)\to (\C^{n+1}\times\C^r,0)$, but when working with map-germs $f\colon (X,S)\to (\C^{n+1},0)$, it is convenient to extend this notion to allow for the source space $X$ to be deformed as well. A possible definition is as follows:

Given a multi-germ $f\colon X\to Y$ between complex spaces, an unfolding of $f$ over $\Delta=(\C^r,0)$ is a multi-germ $\phi\colon \mathcal X\to Y\times \Delta$ such that, taking the commutative diagram 
\[
\begin{tikzcd}[row sep=2.3em,column sep=.6em]
\mathcal X \arrow[rr,"\phi"]	\arrow[dr]&	&  Y\times \Delta\arrow[ld] \\
	& \Delta&
\end{tikzcd}
\]
where $ Y\times \Delta\to \Delta$ is the usual projection from the product, the following conditions are satisfied:
\begin{enumerate}
\item $\mathcal X\to \Delta$ is flat family, with central fibre $\mathcal X_0=\mathcal X\times_\Delta\{0\}$ isomorphic to $X$.
\item The above isomorphism $\mathcal X_0\to X$  produces  a commutative diagram
\[
\begin{tikzcd}[row sep=.8em,column sep=2.3em]
\mathcal X_0 \arrow[rd,"\phi_0"]	\arrow[dd]&\\
&Y \\
X\ar[ru,"f",swap]	&
\end{tikzcd}
\]
\end{enumerate}
In particular, this provides a family of mappings $\mathcal X_t\xrightarrow{\phi_t} Y$ so that $\phi_0$ is equivalent to $f$. 

For our purposes, the crucial property is that  unfoldings of finite multi-germs satisfy 
\[M_k(f)\cong M_k(\phi)\times_\Delta\{0\}.\]
This can be found in more general form in \cite[Proposition 11.6]{Mond-Nuno2020}. Now we are ready to state our main result about triple points in the target.

\begin{thm}\label{thmTriples}Given $f\colon (X,S)\to (\C^{n+1},0)$, assume $f$ to be generically one-to-one and that $M_3^{(i)}=\emptyset$, for $2\leq i\leq r$. If $M_3^{(1)}\neq \emptyset$, then assume $(X,x^{(1)})$ to be Gorenstein. Assume also that $f$ admits an unfolding $\phi$ such that $M_3(\phi)$ and  $M_4(\phi)$ are dimensionally correct. If $M_3(f)$ is dimensionally correct, then it is Cohen-Macaulay. If furthermore $M_4(f)$ is dimensionally correct, then the ideal of $M_3(f)$ satisfies
\begin{align*}\label{thmTriples}
F_2(f_*\cO_{(X,S)})=&F_2^{(1)}\cap \bigcap_{1\leq i<j\leq r}\big(F_1^{(i)}+F_0^{(j)}\big)\cap \big(F_0^{(i)}+F_1^{(j)}\big)\\
				&\cap \bigcap_{1\leq i<j<k\leq r}\big(F_0^{(i)}+F_0^{(j)}+F_0^{(k)}\big).
 \end{align*}
 \begin{proof}
From the isormorphism $M_k(f)\cong M_k(\phi)\times_\Delta\{0\}$ and the fact that $\{0\}$ is regularly embedded in $\Delta$, it follows that, if $M_3(\phi)$ is Cohen-Macaulay, then so is $M_3(f)$, because the hypothesis that $M_3(f)$ is dimensionally correct implies that it is regularly embedded in $M_3(\phi)$. To keep the notation we are used to, instead of replacing $f$ by $\phi$ and showing $M_3(\phi)$ to be Cohen Macaulay, we may just assume $M_3(f)$ and $M_4(f)$ to be dimensionally correct and show $M_3(f)$ to be Cohen-Macaulay.

We proceed by induction on $|S|$ as in the proof of Theorem \ref{thmMultigermDoublePoints}. The case of $r=1$ is precisely Theorem \ref{MondPellikaanTriple} and, for $r\geq 2$, we again use the notation $S'=\{x^{(1)},\dots,x^{(r-1)}\}$ and let \[f'= f_{(X,S')}\quad\text{and}\quad F'_j=F_j(f'_*\cO_{(X,S')}).\] 

The ideal $F_0^{(r)}$ is principal and, since $f^{(r)}$ is generically one-to-one and $M_3^{(r)}(f)=\emptyset$, the ideal $F_1^{(r)}$ is generated by a regular sequence, by Theorem \ref{MondPellikaanDouble}. From the fact that $M_4(f)$ is dimensionally correct, it follows that the generators of $F_0^{(r)}$ and $F_1^{(r)}$ may be chosen so that they form regular sequences in $\cO_{n+1}/F_2'$ and $\cO_{n+1}/F_1'$, respectively. This implies $F_2'F_0^{(r)}=F_2'\cap F_0^{(r)}$ and $F_1'F_1^{(r)}=F_1'\cap F_1^{(r)}$. Then, taking into account that $F_2^{(r)}=\cO_{n+1}$, and exploiting the inclusions $F_0'\subseteq F_1'\subseteq F_2'$ and $F_0^{(r)}\subseteq F_1^{(r)}$, one obtains
\begin{align*}
F_2(f_*\cO_{X,S})	&=F_2' F_0^{(r)}+F_1' F_1^{(r)}+F_0'\\
 				&=F_2' \cap F_0^{(r)}+F_1' \cap F_1^{(r)}+F_0'\\
				&=F'_2 \cap \big(F'_{1}+F^{(r)}_{0}\big)\cap \big(F'_{0}+F^{(r)}_{1}\big).
\end{align*}
The same inclusions also imply that
\[\Big(F'_2 \cap \big(F'_{1}+F^{(r)}_{0}\big)\Big)+ F'_{0}+F^{(r)}_{1}=\Big(F'_2 \cap \big(F'_{1}+F^{(r)}_{0}\big)\Big)+ F^{(r)}_{1} =F_1'+F_1^{(r)}.\]
This gives a short exact sequence
\[
0\to\frac{\cO_{n+1}}{ F_2(f_*\cO_{X,S})  }\to\frac{\cO_{n+1}}{ F'_2 \cap\big(F'_{1}+F^{(r)}_{0}\big) }\oplus \frac{\cO_{n+1}}{  F'_{0}+F^{(r)}_{1} }\to
\frac{\cO_{n+1}}{  F'_{1}+F^{(r)}_{1}  }\to 0.
\]
The hypotheses that $M_3(f)$ and $M_4(f)$ are dimensionally correct imply that \[\depth \frac{\cO_{n+1}}{  F'_{0}+F^{(r)}_{1} }=n-2\quad\text{and}\quad\depth \frac{\cO_{n+1}}{  F'_{1}+F^{(r)}_{1} }=n-3.\] Therefore, by Item (\ref{lemDepthLemmaItem3}) of Lemma \ref{lemDepthLemmaItem3}, it suffices to show that
\[\depth \frac{\cO_{n+1}}{ F'_2 \cap \big(F'_{1}+F^{(r)}_{0}\big)} \geq n-2.\]
Since $F_1'\subseteq F_2'$, this last ring fits into the short exact sequence
\[0\to\frac{\cO_{n+1}}{   F'_2 \cap \big(F'_{1}+F^{(r)}_{0}\big)  }\to\frac{\cO_{n+1}}{ F'_2 }\oplus \frac{\cO_{n+1}}{  F'_{1}+F^{(r)}_{0} }\to \frac{\cO_{n+1}}{  F'_2 +F^{(r)}_{0} }\to 0.\]
That the depth is the desired one follows again by Item (\ref{lemDepthLemmaItem3}) of Lemma \ref{lemDepthLemmaItem3}, taking into account that the dimensionality hypotheses imply \[\depth \frac{\cO_{n+1}}{ F'_2 }=\depth \frac{\cO_{n+1}}{  F'_{1}+F^{(r)}_{0} }=n-2\quad\text{and}\quad\depth  \frac{\cO_{n+1}}{  F'_2 +F^{(r)}_{0} }=n-3.\]
This completes the proof of the statement that $M_3(f)$ is Cohen-Macaulay.

Now assume $M_3(f)$ and $M_4(f)$ are dimensionally correct. We have shown before that, under these hypotheses, 
\[F_2(f_*\cO_{X,S})=F'_2 \cap \big(F'_{1}+F^{(r)}_{0}\big)\cap \big(F'_{0}+F^{(r)}_{1}\big).\]
From the fact that $f$ is generically one-to-one, it follows that $F'_{0}=\bigcap_{1\leq i\leq r-1}F_0^{(i)}$. Applying induction and Theorem \ref{thmMultigermDoublePoints}, we obtain the expression
\begin{align*}
F_2(f_*\cO_{(X,S)})=&F_2^{(1)}\cap \bigcap_{1\leq i<j\leq r-1}\big(F_1^{(i)}+F_0^{(j)}\big)\cap \big(F_0^{(i)}+F_1^{(j)}\big)\\
				&\cap \bigcap_{1\leq i<j<k\leq r-1}\big(F_0^{(i)}+F_0^{(j)}+F_0^{(k)}\big)\\
				&\cap\bigg(\Big(\bigcap_{1\leq i\leq r-1}F_1^{(i)}\cap \bigcap_{1\leq i<j\leq r-1}\big(F_0^{(i)}+F_0^{(j)}\big)\Big)+F^{(r)}_{0}\bigg)\\
				&\cap \Big(\big(\bigcap_{1\leq i\leq r-1}F_0^{(i)}\big)+F^{(r)}_{1}\Big).
\end{align*}
Comparing this to the desired expression shows that it suffices to prove the equality
 \begin{multline*}
\Big(\bigcap_{1\leq i\leq r-1}F_1^{(i)}\cap \bigcap_{1\leq i<j\leq r-1}\big(F_0^{(i)}+F_0^{(j)}\big)\Big)+F^{(r)}_{0}=\\
\bigcap_{1\leq i\leq r-1}\big(F_1^{(i)}+F_0^{(r)}\big)\cap \bigcap_{1\leq i<j\leq r-1}\big(F_0^{(i)}+F_0^{(j)}+F_0^{(r)}\big)
\end{multline*}
%
and the equality
\[\big(\bigcap_{1\leq i\leq r-1}F_0^{(i)}\big)+F^{(r)}_{1}=\bigcap_{1\leq i\leq r-1}\big(F_0^{(i)}+F^{(r)}_{1}\big).\]
The conditions that $M_4(f)$ and $M_3(f)$ are dimensionally correct imply that, given $i,j,k$, with $i\neq j$, no irreducible component of $V(F_0^{(i)}+F_0^{(j)}+F_0^{(r)})$ is contained in $M_2^{(k)}=V(F_1^{(k)})$, no irreducible component of $M_2^{(k)}\cap M_1^{(r)}$ is contained in $V(F_0^{(i)}+F_0^{(j)})$, no irreducible $M_2^{(i)}\cap M_1^{(r)}$ is contained in $M_2^{(j)}$ and no irreducible component of $V(F_0^{(k)}+F_0^{(i)}+F_0^{(r)})$ is contained in $V(F_0^{(k)}+F_0^{(j)})$.
The first equality then follows from Lemma \ref{lemSumInterCommute}. The same dimensional correctness conditions imply that no irreducible component of $M_1^{(i)}\cap M_2^{(r)}$ is contained in $M_1^{(j)}$, for $i\neq j$. Then by Lemma \ref{lemSumInterCommute} the second equality holds as well.
\end{proof}
\end{thm}

In the case where $X$ is a complete intersection, the hypotheses in Theorem \ref{thmTriples} can be relaxed substantially:

\begin{prop}\label{propTriples}
Given a generically one-to-one multi-germ $f\colon (X,S)\to (\C^{n+1},0)$, assume $X$ to be a complete intersection and that $M_3^{(i)}=\emptyset$, for $2\leq i\leq r$. If $M_3(f)$ is dimensionally correct, then it is Cohen-Macaulay. If furthermore $M_4(f)$ is dimensionally correct, then $M_3(f)$ can be computed with the same formula as in Theorem \ref{thmTriples}.
\begin{proof}The hypothesis that $X$ is a complete intersection implies that $X$ is Gorenstein. To see that $f$ admits an unfolding whose spaces $M_3$ and $M_4$ are dimensionally correct we use an original idea of Mond and Montaldi \cite{Mond-Montaldi:1994}. Take a regular sequence $h_1,\dots,h_k$ defining $(X,S)$ in its ambient space  $(\C^{n+k},S)$ and consider the mapping \[\varphi\colon (\C^{n+k},S)\to(\C^{n+1}\times\C^{k},0),\]  
given by 
\[x\mapsto(\bar{f}_1(x),\dots,\bar{f}_{n+1}(x),h_1(x),\dots,h_k(x)),\]
where $\bar{f}_i$ is any analytic extension in $\cO_{\C^{k+n},S}$ of the $i$th coordinate function $f_i$ of $f$. The mapping $\varphi$ is finite and, from the fact that $h_1,\dots,h_k$ are a regular sequence defined on the smooth space $(\C^{n+k},S)$, it follows that $h_1,\dots,h_k$ define a flat morphism \cite[Theorem 23.1]{Matsumura:1989}, hence $\varphi$ is an unfolding of $f$. Since $\varphi$ is finite, in particular it is $\sK$-finite and, by a construction due to Mather (see e.g. \cite[Theorem 7.2]{Mond-Nuno2020}), it admits a stable unfolding, 
\[F\colon (\C^{k+n}\times\C^r,S\times\{0\})\to(\C^k\times\C^{n+1}\times\C^r,0).\]

Being an unfolding of $\varphi$, the germ $F$ is an unfolding of $f$ as well. Therefore, it suffices to show that the multiple point spaces $M_3(F)$ and $M_4(F)$ are both dimensionally correct when $F$ is stable. 
Fix a small enough finite representative $F\colon U\to V$, where $U,V$ are open neighbourhoods of $S\times \{0\}$ and $0$ in $\C^{k+n}\times\C^r$ and $\C^k\times\C^{n+1}\times\C^r$, respectively. The set $\Sigma^2$ of points in $U$ of corank $\ge 2$ is closed analytic and has codimension $6$ in $U$. In $F(U)\setminus F(\Sigma^2)$ we consider the stratification by stable types, whose strata are the the isosingular loci $\operatorname{Iso}(F,y)$, with $y\in F(U)\setminus F(\Sigma^2)$. By construction, $\operatorname{Iso}(F,y)$ is the set of points $y'\in F(U)$ such that the multi-germ of $F$ at $y'$ is $\sA$-equivalent to the multi-germ of $F$ at $y'$. It follows that $\operatorname{Iso}(F,y)$ is a submanifold of $V$ of codimension equal to the $\sK_e$-codimension of the multi-germ of $F$ at $y$ (see \cite[Theorem 7.4]{Mond-Nuno2020}). Moreover, the $\sA$-class of $F$ at $y$ is determined by the contact algebra $Q_y(F)$ and since $y\notin F(\Sigma^2)$, it must be of the form
\[
Q_y(F)\cong \frac{\C\{x\}}{(x^{a_1})}\oplus\dots\oplus \frac{\C\{x\}}{(x^{a_\ell})},
\]
for some $a_1,\dots,a_\ell\ge 1$. The $\sK_e$-codimension is $(2a_1-1)+\dots+(2a_\ell-1)$ which is the codimension of the stratum $\operatorname{Iso}(F,y)$ in $V$. On the other hand, $y\in M_r (F)$ if and only if $\dim_\C Q_y(F)\ge r$ and we also have 
$\dim_\C Q_y(F)=a_1+\dots+a_\ell$. In conclusion, there exists only one stratum of codimension $r$ in $M_r(F)$ corresponding to the transverse $r$-multiple points and all other strata have codimension $>r$. This shows that $M_r(F)$ is dimensionally correct for $r\le 7$.
\end{proof}
\end{prop}

\begin{ex}\label{exBigerm1}We are going to compute $M_{3}(f)$, for the bi-germ \[f: U^{(1)}\sqcup U^{(2)} \to (\C^{4},0)\]  with branches
\smallskip
\begin{itemize}
\item[] $f^{(1)}\colon U^{(1)}=(\C^{3},0) \to (\C^{4},0), \ 
	(t,x,y)\mapsto (t,x,y^{3}+ty,xy+y^{5})$
	\medskip
\item[] $f^{(2)}\colon  U^{(2)}=(\C^{3},0) \to (\C^{4},0), \ 
	(\tilde{t},\tilde{x},\tilde{y})\mapsto (\tilde{t},\tilde{x},\tilde{y},\tilde{t})$.
\end{itemize}
	The presentation matrices of $f^{(1)}$ and $f^{(2)}$ are 
	\[
	\lambda^{(1)}=
	\begin{bmatrix}
		-Z-TY & Y^{2} & XY-TZ \\
		X+T^{2} & -Z-2TY & Y^{2} \\
		Y  &  X+T^{2} & -Z-TY
	\end{bmatrix},\quad
	\lambda^{(2)}=
	\begin{bmatrix}
		-Z+T
	\end{bmatrix}.
	\]
	
	In order to use the formula in Theorem \ref{thmTriples}, we must check first $M_{4}(f)$ to be dimensionally correct. The presentation matrix of $f$ is the block diagonal matrix
	\[
	\lambda=
	\begin{bmatrix}
		\lambda^{(1)} & 0 \\
		0 &  \lambda^{(2)}
	\end{bmatrix}.
	\]
We obtain $F_{3}=\langle X,Y,Z,T\rangle$, hence $M_4(f)$ has dimension zero, as desired. Since $\lambda^{(2)}$ is a $1\times 1$ matrix, we see that $F_1^{(2)}=\cO_4$, hence  $M_3(f)$ is defined by the ideal
	\[F_2=F_{2}^{(1)} \cap (F_{1}^{(1)}+  F_{0}^{(2)}).\]
It turns out that this is the primary decomposition of $F_2$, expressing it as the intersection of the two prime ideals
\smallskip
\begin{itemize}
\item[] 	$F^{(1)}_{2}=\langle X+T^2,Z,Y\rangle$,
\medskip
\item[] 	$\begin{aligned}[t]F^{(1)}_{1}+ F^{(2)}_{0}=&\langle T-Z,
		X^2+TY+2XT^2+2TY^2+T^4,\\
		&XT+XTY+T^3+Y^3+T^3Y,\\
		&T^2-XY^2+3T^2Y+T^2Y^2
		\rangle .\end{aligned}$
\end{itemize}

We see that $M_3(f)$ is a curve (i.e., it has dimension one) and it is not difficult to check that it is reduced. So, it makes sense to compute its Milnor number $\mu(M_3(f))$ in the sense of \cite[Definition 1.1.1]{Buchweitz-Greuel}. To do it, we use the Milnor formula $\mu=2\delta-r+1$, where $\delta$ is the delta invariant and $r$ is the number of branches (see \cite[Proposition 1.2.1]{Buchweitz-Greuel}). 

We split $M_3(f)$ as a union of two curves $M_3(f)=C\cup D$, where 
\[
C=V(F^{(1)}_{2}),\quad D= V(F^{(1)}_{1}+ F^{(2)}_{0}).
\]
Since $C$ is a regular branch, it has $\delta=0$. The other branch $D$ is isomorphic to its projection $\tilde D$ onto the 3-space $X,Y,T$. To compute the invariants of $\tilde D$, we polar multiplicities, following \cite[Corollary 3.2]{Nuno-Tomazella}:
\[
\mu=m_1-m_0+1,
\] 
where $m_0$ is the multiplicity and $m_1$ is the first polar multiplicity, defined as the Milnor number of a generic linear function $p\colon \C^3\to\C$, $p(X,Y,T)=b_1 X+b_2 Y+b_3 T$, $b_i\in\C$, restricted to the curve. When the curve is smoothable, this is equal to the number of critical points of $p$ on the smoothing of the curve. With the aid of \textsc{Singular} \cite{DGPS}, we see that $\tilde D$ admits a determinantal representation as $\tilde D=V(I_2(M))$, where $I_2(M)$ is the ideal generated by the $2\times 2$ minors of the matrix
\[
M=
\left(
\begin{array}{cc}
Y^2 &   -T-YT\\
-T-2YT & X+T^2\\
X+T^2 & Y 
\end{array}
\right).
\]
A smoothing $\tilde D_s$ of $\tilde D$ can be obtained as $\tilde D_s=V(I_2(M+sA))$ where $A=(a_{ij})$ is a generic $3\times 2$-matrix with entries $a_{ij}$ in $\C$ and $s\in\C\setminus\{0\}$ is small enough (see \cite[Lemma 3.1]{Nuno-Orefice-Tomazella}). Again with \textsc{Singular}, we compute the number of critical points of $p$ on $\tilde D_s$, which gives $m_1=6$, and the multiplicity $m_0=3$. This gives
\[
\mu(D)=\mu(\tilde D)=m_1-m_0+1=4.
\]
Since $D$ is irreducible, this implies $\delta(D)=2$. 

Now $\delta(M_3(f))$ follows from a formula due to Hironaka (see \cite[Lemma 1.2.2]{Buchweitz-Greuel}):
\[
\delta(M_3(f))=\delta(C)+\delta(D)+C\cdot D
\]
where $C\cdot D=\dim_\C \O_4/(F^{(1)}_{2}+(F^{(1)}_{1}+ F^{(2)}_{0}))$ is the intersection number of the two branches. Again with \textsc{Singular} we get
$C\cdot D=1$. Thus,
\[
\delta(M_3(f))=0+2+1=3
\]
and hence
\[
\mu(M_3(f))=2\delta(M_3(f))-r+1=5.
\]

\end{ex}
\section{Triple and quadruple points in the source}

Given a finite multi-germ $f\colon (X,S)\to (\C^{n+1},0)$, the multiple point spaces $M_k(f)$ of the previous section lie on the target of $f$, that is, on $(\C^{n+1},0)$. There are however other multiple point spaces that one can define, and the ones this section is devoted are subspace of the source of $f$, instead of its target. For simplicity, we consider only the case where the source is smooth. Given a multi-germ
 \[f\colon (\C^n,S)\to (\C^{n+1},0),\]
 the idea is as follows: Start by defining a subspace $D^2(f)\subseteq (\C^n,S)\times (\C^n,S)$ which, naively, should be the closure of the set of pairs $(x,x')$ with $x\neq x'$ and $f(x)=f(x')$ (the proper definition is given below). If $f$ is finite, then the projection $(\C^n,S)\times (\C^n,S)\to (\C^n,S)$ on the first coordinate restricts to a finite mapping
\[
\begin{tikzcd}
D^{2}(f)\arrow{d}{p} \\
(\C^n,S)
\end{tikzcd}
\]
The \emph{source double point space} $D(f)=D^2_1(f)$ is defined as the image of $p$, with the analytic structure given by the $0$th Fitting ideal, that is,
\[D^2_1(f)=M_1(p).\]
This makes sense since, given a pair $(x,x')$, with $x\neq x'$ and $f(x)=f(x')$, the point $x\in \C^n$ should be considered a double point of $f$.

 The same idea tells us that $M_2(p)\subseteq (\C^n,S)$ should be regarded as a triple point space, because two pairs $(x,x')$ and $(\tilde x, \tilde x')$, with $x'\neq x$, $f(x')=f(x)$, $\tilde x'\neq \tilde x$, $f(\tilde x')=f(\tilde x)$, $(x,x')\neq (\tilde x, \tilde x')$ and $x=\tilde x$, carry the same information than a tuple $(x,x',x'')$, with $x\neq x'\neq x''\neq x$ and $f(x)=f(x')=f(x'')$, by setting $x''=\tilde x'$, and then $x$ should be a triple point of $f$. By the same token, points in $M_3(p)$ should be regarded as quadruple points in the source. Before making this into a formal definition, we must define the space $D^2(f)$. 
 
 Consider the ideal $I_\Delta$ defining the diagonal $(\Delta\C^{n+1},0)$ in $(\C^{n+1},0)\times (\C^{n+1},0)$. The fibered product $(\C^n,S)\times _{(\C^{n+1},0)} (\C^n,S)$ is the zero locus of the pullback $(f\times f)^*I_\Delta$. Given  $x^{(i)},x^{(j)}\in S$, we may take coordinates $x$ and $x'$ of $\C^n$, centered at $x^{(i)},x^{(j)}$, so that 
 \[((f\times f)^*I_\Delta)_{(x^{(i)},x^{(j)})}=\big\langle f^{(i)}_1(x')-f^{(j)}_1(x),\dots, f^{(i)}_{n+1}(x')-f^{(j)}_{n+1}(x)\big\rangle.\]
 
 Away from the diagonal, $(\C^n,S)\times _{(\C^{n+1},0)} (\C^n,S)$ corresponds with the idea of what a double point should be, but around a point $(x^{(i)},x^{(i)})$, it contains all points $(x,x)\in \Delta \C^n$, which, a priori, should not be regarded as triple points. To remedy this, we must add extra conditions on the diagonal. If we think of $f^{(i)}(x')-f^{(i)}(x)$ and $x'-x$, the fact that $f^{(i)}(x')-f^{(i)}(x)$ vanishes on $(\Delta \C^n,(x^{(i)},x^{(i)}))$ translates, by Hilbert Nullstellensatz, into  the fact that there is an $(n+1)\times n$ matrix $\alpha$, with entries in $\cO_{2n}$, such that
\[f^{(i)}(x')-f^{(i)}(x)=\alpha^{(i)}(x,x')\cdot(x'-x).\]
The \emph{double point space} $D^2(f)$ is the disjoint union of the germs \[\left(D^2(f),{(x^{(i)},x^{(j)})}\right)=V(I^2_{ij}),\] where 
\[I_{ij}=
\begin{cases}
((f\times f)^*I_\Delta)_{(x^{(i)},x^{(j)})}&\text{ if $i\neq j$}\\
((f\times f)^*I_\Delta)_{(x^{(i)},x^{(i)})}+\big\langle \text{$n\times n$ minors of $\alpha$},\big\rangle&\text{ if $i=j$}
\end{cases}
\]

The space $D^2(f)$ has dimension at least $n-1$ and it is Cohen-Macaulay when it has dimension exactly $n-1$ \cite[Section 9.4]{Mond-Nuno2020}. If $f$ has corank one at $x^{(i)}$, then, after a change of coordinates in source and target, $f^{(i)}$ takes the form of the germ at the origin of a mapping of the form 
\[(x_1,\dots,x_{n-1},y)\mapsto ( x_1,\dots,x_{n-1},f_n(x,y),f_{n+1}(x,y)).\] Then,  $\left(D^2(f),{(x^{(i)},x^{(i)})}\right)$ is isomorphic \cite[Section 9.5]{Mond-Nuno2020} to germ of
\[\{(x,y,y')\in \C^n\times \C\mid \frac{f_n(x,y')-f_n(x,y)}{y'-y}=\frac{f_{n+1}(x,y')-f_{n+1}(x,y)}{y'-y}=0\}.\]
From this, it follows that $\left(D^2(f),{(x^{(i)},x^{(i)})}\right)$ is a complete intersection whenever it has dimension $n-1$.

\begin{definition}
For a finite map-germ $f\colon (\C^n,S)\to (\C^{n+1},0)$, we define the \emph{source triple point space} of $f$ as
\[D^2_2(f)=M_2(p),\]
where $p\colon D^2(f)\to (\C^n,S)$ is as above. Similarly, we define the \emph{source quadruple point space} of $f$ as
\[D^2_3(f)=M_3(p).\]

\end{definition}

\begin{ex} We are going to compute the source triple point space $D^2_2(f)=M_2(p)$ and source quadruple space $D^2_3(f)=M_3(p)$, for the map-germ $f\colon U^{(1)}\sqcup U^{(2)}\to(\C^4,0)$  of Example \ref{exBigerm1}. We also show how the explicit formula for $F_1(f_*\cO_{(X,S)})$ in Theorem \ref{thmMultigermDoublePoints} helps us compute the Milnor number of $D^2_2(f)$ and how Theorem \ref{thmTriples} allows us to compute the number of quadruple points in a stabilization of the multi-germ in question.
Recall that there is a diagram
	\[
	\begin{tikzcd}
		D^{2}(f)\arrow{d}{p} \\
		U^{(1)}\sqcup U^{(2)} \arrow{r}{f} & (\C^{4},0)
	\end{tikzcd}
	\]

The double space of $f$ decomposes as $D^{2}(f)= X_{1} \sqcup X_{2} \sqcup X_{3}$, with $X_1=D^{2}(f^{(1)})$, $X_2=U^{(1)}\times_{(\C^4,0)} U^{(2)}$ and $X_3=U^{(2)}\times_{(\C^4,0)} U^{(1)}$. These spaces are
\begin{itemize}
\item[]	
	$\begin{aligned}
	X_{1}= \Big\{ (t,x,y,y') \mid &\frac{xy'+y'^{5}-xy-y^{5}}{y'-y}= \frac{y'^{3}+ty'-y^{3}-ty}{y'-y}=0\Big\} ,
	\end{aligned}$
\medskip
\item[] 
	$\begin{aligned}
	X_{2}&=\left\{ (t,x,y,\tilde{t},\tilde{x},\tilde{y}) \ | \ \tilde{t}=t, \tilde{x}=x, \tilde{y}=y^{3}+ty,\tilde{t}=xy+y^{5} \right\}\\
		&\cong \left\{ (t,x,y,\tilde{y}) \ | \  y^{3}+ty-\tilde{y}=xy+y^{5}-t=0\right\},
		\end{aligned}$
\medskip
\item[]
	$\begin{aligned}
	X_{3}&=\left\{ (\tilde{t},\tilde{x},\tilde{y},t,x,y) \  | \ \tilde{t}=t, \tilde{x}=x, \tilde{y}=y^{3}+ty,\tilde{t}=xy+y^{5}\right\}\\
		&\cong \left\{ (\tilde{y},t,x,y) \  | \  y^{3}+ty-\tilde{y}=xy+y^{5}-t=0\right\}.
	\end{aligned}$

\end{itemize}

 The projection $p\colon D^2(f)\to U^{(1)}\sqcup U^{(2)}$ decomposes as \[p=\pi \sqcup \tilde{\pi}\colon (X_1\sqcup X_2)\sqcup X_3\to U^{(1)}\sqcup U^{(2)},\] where $\tilde{\pi}: X_{3} \to U^{(2)}$ is a mono-germ and \(\pi\) is a bi-germ with two branches $\pi^{(1)}:X_{1} \to U^{(1)}$ and $\pi^{(2)}:X_{2} \to U^{(1)}$. Therefore, the triple and quadruple points in the source of $f$ decompose as
 \[D^2_2(f)=M_2(\pi)\sqcup M_2(\tilde\pi)\quad\text{and}\quad D^2_3(f)=M_3(\pi)\sqcup M_3(\tilde\pi).\]
 We start with the computations over $U^{(1)}$, by computing some presentation matrices of $\pi^{(1)}_*\cO_{X_1}$ and  $\pi^{(2)}_*\cO_{X_2}$, namely
		\[\xi^{(1)}=
	\begin{bmatrix}
		y^{4}+y^{2}t+t^{2}+x &y^{5}+2y^{3}t+yt^{2}\\
		 y^{3}+yt  & -t^{2}-x
	\end{bmatrix},\quad 
	\xi^{(2)}=
	\begin{bmatrix}
		y^{5}+xy-t
	\end{bmatrix}. \]
Since $M_{4}(p)$ on $U^{(1)}$ is empty and $X_1\sqcup X_2$ is a complete intersection, Proposition \ref{propTriples} allows us to use the formula in Theorem \ref{thmTriples} for $\pi$, which in this case reads	\[M_{3}(\pi)=V\big(F_{1}^{(1)} + F_{0}^{(2)}\big)=V(x,t,y^3).\] 	As \( \dim M_{3}(\pi) =0 \), we may apply  the formula in Theorem \ref{thmMultigermDoublePoints}, which gives
	\[
	M_{2}(\pi)=V\Big(F_{1}^{(1)} \cap \big(F_{0}^{(1)}+F_{0}^{(2)}\big)\Big),
	\]
with
\smallskip
\begin{itemize}
\item[] $F_{1}^{(1)}=  \langle y,x+t^{2}\rangle \cap \langle x+t^{2},t+y^{2}\rangle$,
\medskip
\item[] 	$\begin{aligned}[t] F_{0}^{(1)}+F_{0}^{(2)}=\langle& t-xy-y^5,\\ &x^2+2xt^2+3yt^2-2xy^2t+y^3t+t^4+2y^2t^3+4y^4t^2
		\rangle.
	\end{aligned}$
\end{itemize}

Before moving on to the calculations  over $U^{(2)}$, we compute $\mu(M_2(\pi))$.  We know $M_2(\pi)$ decomposes as the union of $V(F_1^{(1)})$ and $F_{0}^{(1)}+F_{0}^{(2)}$ and that, in turn, $V(F_1^{(1)})$ decomposes as the union of two regular curves, which we call $Z_1$ and $Z_2$. In turn, the expression of $F_{0}^{(1)}+F_{0}^{(2)}$ allows us to eliminate $t$, turning $V(F_{0}^{(1)}+F_{0}^{(2)})$ into a plane curve with Milnor number $7$. This, combined with Hironaka's $\mu=2\delta-r+1$ formula, prevents this curve to be irreducible and, looking at the degree two part, we conclude that it consists of $2$ regular curves, called $Z_3$ and $Z_4$, with $\delta(Z_3\cup Z_4)=4$. Computing the intersection numbers $Z_1\cdot (Z_2\cup Z_3\cup Z_4)=2$ and $Z_2\cdot (Z_3\cup Z_4)=2$, we conclude
\begin{align*}
\delta(M_2(\pi))	=&\delta (Z_1)+\delta (Z_2)+\delta(Z_3\cup Z_4)\\
			&+Z_1\cdot (Z_2\cup Z_3\cup Z_4)+Z_2\cdot (Z_3\cup Z_4)\\
			=&0+0+4+2+2=8
\end{align*}
Now from Hironaka's formula we conclude that
\[\mu(M_2(\pi))=2\cdot 8-4+1=13.\]

Now we proceed to the calculations over $U^{(2)}$. A presentation matrix of $\tilde{\pi}_*\cO_{X_3}$ is 
		\[ 
	\tilde{\xi}=
	\begin{bmatrix}
		-\tilde{y} & 0 & t & 0 & t^{2} \\
		 t & -\tilde{y} & -x  & t & -xt\\
		 0 & t & -\tilde{y}  & -x & t\\
		 1 & 0  & t& -\tilde{y} & -x\\
		 0 & 1 & 0 & t &-\tilde{y}
	\end{bmatrix}
	\]
The source triple space $D^2_2(f)\cap U^{(2)}=M_{2}(\tilde \pi)=V(\tilde F_1)$ is reduced and irreducible, with defining ideal
\smallskip
\begin{itemize}
\item[]	
$\begin{aligned}[t]\tilde{F_{1}}=& \langle t^2+3t^2\tilde{y}-x\tilde{y}^2+t^2\tilde{y}^2,
		tx+t^3+tx\tilde{y}+\tilde{y}^3+t^3\tilde{y},\\
		&x^2+t\tilde{y}+2t^2x+2t\tilde{y}^2+t^4
		 \rangle
\end{aligned}$
\end{itemize}
We observe that $M_{2}(\tilde \pi)$ coincides with the curve $\tilde D$ in Example \ref{exBigerm1}, hence $\mu(M_{2}(\tilde \pi))=4$.

The source quadruple space is  \[D^2_3(f)\cap U^{(2)}=M_{3}(\tilde \pi)=V(\tilde{F_{2}})= V( x, t, \tilde{y} ).\]

We may use these calculations to compute the number of quadruple points found in a stabilization of $f$, namely
\[Q(f)=\frac{\dim_\C\cO_{D^2_3(\pi)}+\dim_\C\cO_{D^2_3(\tilde \pi)}}{4}=\frac{3+1}{4}=1.\] 
To understand that this  is the number of quadruple points in a stabilization, one first observes that, if $(f_t,t)$ is a stabilization of $f=f_0$, then Proposition \ref{propTriples} applies to the mapping $D^2(f_t,t)\to (U^{(1)}\sqcup U^{(2)})\times(\C,0)$ and, therefore, the space $D^2_3(f_t,t)$ is a reduced Cohen Macaulay curve. The classification of stable multi-germs $(\C^3,S)\to(\C^4,0)$ implies that, for the stable mapping $f_\epsilon,\epsilon \neq 0$, the space $D^2_3(f_\epsilon)$ consist of four reduced points for each quadruple point. Now  the claim follows by the principle of conservation of number.

As a finally observation, it seems natural to wonder how the triple and quadruple point spaces $D^2_2(f)$ and $D^2_3(f)$, which are subspaces of $(X,S)$, compare to the preimages by $f$ of $M_3(f)$ and $M_4(f)$, which are subspaces of $(\C^{4},0)$ and were computed in Example \ref{exBigerm1}. A simple computation shows, for the multi-germ $f$ in question, the equalities $D^{2}_2(f)=f^{-1}(M_{3}(f))$ and $D^2_{3}(f)=f^{-1}(M_{4}(f))$ hold.

\end{ex}

\bibliography{CMTargetDoublePoints} 
\bibliographystyle{plain}	
%
\end{document}